\documentclass[12pt]{article}

\begin{document}

\title{Examples of centralizers in the Artin braid groups}         
\author{Nikolai V. Ivanov}        
\date{}          

\maketitle

\section{Introduction}

In the first version of the preprint \cite{FG-M}, N. Franco and J. Gonz\'alez-Meneses made
the following conjecture (discussed also in the second version \cite{FG-M}):\\

{\it The normalizer of an element of the Artin braid group $B_n$ on $n$ strings
is generated by no more than $n-1$ elements.}\\

{\it The normalizer} of an element of a group is defined as the subgroup of all
elements commuting with it; this subgroup is often also called {\it the centralizer}
of the given element. In the second version \cite{FG-M} the authors
switched to the latter terminology, and we will also use the term ``centralizer''.

The conjecture was supported by a new algorithm for finding generators
of a centralizer, suggested in \cite{FG-M}, and by extensive computations
based on this algorithm.

The goal of this paper is to present simple examples of elements
of $B_n$ for which the centralizer cannot be generated by a less than
quadratic in $n$ number of generators. In particular, the above conjecture
is disproved. The key insight is that the centralizer may be closer to
a pure braid group than to a braid group (as one may think initially), 
and that a pure braid group requires the number of generators quadratic
in the number of strings.

Our methods are based on Thurston's theory of the
surface diffeomorphisms and on the well known relation between the
Artin braid groups and mapping class groups, in contrast with the 
purely algebraic methods of \cite{FG-M}. In fact, 
the effectiveness of  Thurston's theory for studying the centralizers
in the context of the mapping class groups was established 
long ago by J. McCarthy \cite{M} and the author \cite{I1}, \cite{I2}.
Centralizers, considered from the point of view of Thurston's
theory, also play a key role in the recent work \cite{IIM} about braid
groups of surfaces of higher genus of E. Irmak, J. McCarthy and the author.

One should also mention a nearly forgotten paper of G. Burde \cite{B}
in which he gave, in particular, a description of the centralizers
of the elements of the pure braid group in the whole Artin braid
group.  The methods of Burde are fairly geometrical; Thurston's theory 
clarifies the picture further. It would be interesting to compare in detail
the methods of Burde and the methods based on Thurston's theory.

The present work was reported in the author's
talk \cite{I3}. After learning about these results, J. Gonz\'alez-Meneses 
jointly with B. Wiest \cite{GM-W} applied Thurston's theory
to give a rather complete description of the centralizers in
the Artin braid group and found a conjecturally exact upper bound for
the minimal number of generators of a centralizer. This bound is exact
for the examples discussed in the present paper.

J. Gonz\'alez-Meneses and B. Wiest \cite{GM-W}, while using 
Thurston's theory, present many of their arguments in the
framework of braids. In contrast, in the present note we work 
entirely in the framework of the mapping class groups. One may
hope that this point of view may clarify some of the more
obscure parts of \cite{GM-W} (the authors themselves point out
that some proofs seem to be more complicated than necessary). 
 
\section{Preliminaries}

Let $B_n$ be the Artin braid group on $n$ strings. The group $B_n$ 
can be realized as the group of isotopy classes of
diffeomorphisms of the disc $D^2$ pointwise fixed on the boundary
$\partial D^2$ and preserving a set of $n$ distinguished points
in ${\rm int}\, D^2$, called {\em punctures}. This result is, essentially,
due to Artin himself \cite{Art}. It is usually assumed that $D^2$ is the 
standard unit disc in ${\bf R}^2$ and the punctures are lying on the real 
line; this allows a canonical identification of the group $B_n$ defined
in terms of Artin generators and relations with the above group of isotopy
classes. 

If we remove the condition that the diffeomorphisms of $D^2$ are pointwise fixed
on the boundary, the resulting group of isotopy classes is known as
the mapping class group of the disc with $n$ punctures and will be denoted
by $M_n$.  Obviously, there is a canonical forgetting homomorphism
\[B_n \rightarrow M_n,\]
where $B_n$ is realized, as above, as the group of isotopy classes.
The following facts are well known.
First, this homomorphism is surjective. Its kernel is equal to 
the center of $B_n$ and is an infinite
cyclic group generated by the Dehn twist about the boundary $\partial
D^2$.

This relation between $B_n$ and $M_n$ shows that in order to find elements 
with the centralizers generated by, say, no less than a quadratic in $n$ number 
of elements it is sufficient to find elements with the same property in $M_n$. 

One may also the consider the {\em pure} braid group $PB_n$ and the corresponding
mapping class group $PM_n$; in order to define them in terms of the isotopy
classes, one needs only to add the condition that the considered
diffeomorphisms preserve each puncture individually. Similarly to the
above, there is a surjective homomorphism $PB_n \rightarrow PM_n$ 
having central kernel isomorphic to $\bf Z$. A crucial
fact about pure braid groups is the following theorem of V. I. Arnold \cite{Arn}.\\

{\em The first homology group $H_1 (PB_n, {\bf R})$ is isomorphic to 
${\bf R}^{n(n-1)/2}$.}\\

This theorem immediately implies that the group $PB_n$ cannot be generated
by less than $n(n-1)/2$ elements, and its quotient group $PM_n$ by less
than $\frac{n(n-1)}{2}-1$ elements.

\section{Examples} 

Let us assume for simplicity that $n$ is divisible by $3$, say $n=3m$. 
Let $D_1, \ldots, D_m$ be $m$ disjoint discs in the interior of $D^2$ 
such that every disc $D_i$, $1\leq i \leq m$,
contains exactly $3$ punctures in its interior. 
Pick up for each $i$, $1\leq i \leq m$,
a pseudo-Anosov isotopy class $f_i$ of diffeomorphisms $D_i \rightarrow D_i$.
We can choose such $f_i$'s that the dilatation coefficients (assumed to be 
$>1$ of all these pseudo-Anosov classes are different. Represent each of 
these pseudo-Anosov classes by a diffeomorphism
$D_i \rightarrow D_i$ pointwise fixed on the boundary $\partial D^2$, and
extend these diffeomorphisms by the identity to a diffeomorphism $F$ of $D^2$.
Let $f$ be the isotopy class of $F$.

The Thurston's normal form of $f$ is obvious: $f$ is reducible; boundaries
$\partial D_i$ form the canonical reduction system of $f$; the isotopy classes 
$f_i$ and the identity class on the complement of the interiors of discs $D_i$
are the canonical ``pieces'' of $f$. 

If $g\in M_n$, then the Thurston's normal
form of $gfg^{-1}$ results from applying (a representative of) $g$ to the
Thurston's normal form of $f$. In particular, if $g$ belongs to the centralizer of $f$,
i.~e., if $gfg^{-1} = f$, then $g$ preserves the (isotopy class of) the union of
the boundaries $\partial D_i$. Moreover, since the pseudo-Anosov classes
$f_i$ have different dilatation coefficients, $g$ must preserve every
boundary $\partial D_i$ and every disc $D_i$ individually. The
structure of the centralizer is now clear: it consists of isotopy classes of
diffeomorphisms $G$ preserving discs $D_i$ and subject only to the condition
that the isotopy classes of the restrictions on discs $D_i$ commute with
the classes $f_i$. This condition means that these restrictions belong
to the centralizers of $f_i$'s, and one can use the results of J. McCarthy
\cite{M} to describe these centralizers. 

More important for us is the
fact that there is no condition at all on the action of $G$ on the complement
of the interiors of $D_i$'s, except that all boundary components are preserved.
By restricting $G$ to this complement and collapsing  boundaries $\partial D_i$
to punctures, we get a surjective homomorphism of the centralizer of $f$ to 
the group $PM_m$ of the isotopy classes of diffeomorphisms of a disc
fixing each of $m$ punctures in the interior. 

Using the above mentioned corollary of the Arnold's theorem, we see that
$PM_m$ cannot be generated by less than $\frac{m(m-1)}{2}-1$ elements. Since
there is a surjective homomorphism from the centralizer of $f$
to $PM_m$, this centralizer also cannot be generated by less than
$\frac{m(m-1)}{2}-1$ elements. In fact, an easy argument shows that
for any lift of $f$ to $PB_n$ its centralizer admits a surjective
homomorphism to $PB_m$, and hence cannot be generated by less than
$m(m-1)/2$ elements.

In any case, since $n=3m$, our elements $f$ require
a quadratic in $n$ number of elements to generate its centralizer.
This, clearly, disproves the conjecture stated in Introduction.

\section{Remarks}

\paragraph{1.} Notice that the free abelian group generated by $m$ Dehn twists
about boundaries $\partial D_i$ is obviously contained in the center
of the centralizers of our examples $f$. This implies that the centralizer
of a lift of $f$ to $PB_n$ cannot be generated by less than $\frac{m(m-1)}{2}+m=
\frac{m(m+1)}{2}$ elements.

\paragraph{2.} The above examples admit a lot of variations.
H. Hamidi-Tehrani suggested (right after the talk \cite{I3}) 
a modification of these examples: instead of pseudo-Anosov 
classes with different dilatations one can use different
powers of Dehn twists about the boundaries $\partial D_i$.
After this, it is only natural to notice that there is no
need to have $3$ punctures inside the discs ($3$ punctures
were needed to have pseudo-Anosov classes on $D_i$'s), and therefore 
one can have exactly $2$ punctures inside every $D_i$. Also, one may
leave (exactly) one puncture outside all these discs, if the number $n$ of
punctures is odd. This leads
to examples of elements requiring at least $\frac{m(m-1)}{2}+m=\frac{m(m+1)}{2}$ elements
to generate their centralizer for $n=2m$ and $\frac{m(m+1)}{2}+m=\frac{m(m+3)}{2}$ elements
for $n=2m+1$. According to \cite{GM-W}, the first specific examples of this type (with
$2$ punctures inside the discs) were suggested by S.~J. Lee. It is proved in \cite{GM-W} that
such examples are the worst possible in terms of the required number of
generators of the centralizer.

\paragraph{3.}Another illustration of the power of the geometric approach to this 
circle of questions is provided by some specific elements of Artin braid 
groups for which a computation of the centralizers was done by methods of
the combinatorial group theory by G. G. Gurzo \cite{Gu1}, \cite{Gu2}.
The proofs are based on long and difficult computations.
All elements considered by G. G. Gurzo have simple images in $M_n$,
and Thurston's theory very easily leads to a description of their
centralizers. We leave this as an exercise for an interested reader.
Many of the specific elements treated by Gurzo can be
treated also by the methods of G. Burde \cite{B}; 
strangely enough, Burde's work is cited by Gurzo, but is not used.

\noindent
{\sc Michigan State University\\
Department of Mathematics\\
Wells Hall\\
East Lansing, MI 48824-1027\\

\noindent
E-mail:} ivanov@math.msu.edu


\begin{thebibliography}{XXXX}

\bibliography{References}


\bibitem[Arn]{Arn} V. I. Arnold, The cohomology ring of the group of dyed braids,
{\it Mat. Zametki}, V. 5 (1969), 227-231.

\bibitem[Art]{Art} E. Artin, Theory of braids, Ann. of Math., V. 49 (1948), 101-126.

\bibitem[B]{B} G. Burde, \"Uber Normalisatoren der Zopfgruppe,
Abh. Math. Sem. Univ. Hamburg, v. 27 (1964), 235-254.

\bibitem[F-GM]{FG-M} N. Franco, J. Gonz\'alez-Meneses, Computation
of centralizers in braid groups and Garside groups, Preprint, 
arXiv:math.GT/0201243. To appear in {\it Revista matem\'atica Iberoamericana}.

\bibitem[GM-W]{GM-W} J. Gonz\'alez-Meneses, B. Wiest, On the structure
of the centralizer of a braid, Preprint, arXiv:math.GT/0305156.

\bibitem[Gu1]{Gu1} G. G. Gurzo,  Systems of generators for normalizers of certain
elements of the braid group, Izvestia AN SSSR, v. 48, No. 3 (1984), 476-519.
{\bf English translation:} Mathematics of the USSR-Izvestia, v. 24, No. 3 (1985)
439-478.

\bibitem[Gu2]{Gu2} G. G. Gurzo, Systems of generators for centralizers of rigid elements
of the braid group, Izvestia AN SSSR, v. 51, No. 5 (1987), 915-935.
{\bf English translation:} Mathematics of the USSR-Izvestia, v. 24, No. 3 (1985)
223-244.

\bibitem[IIM]{IIM} E. Irmak, N. V. Ivanov, J. D. MacCarthy, Automorphisms 
of surface braid groups, Preprint, arXiv:math.GT/0306069.

\bibitem[I1]{I1} N. V. Ivanov, Automorphisms of Teichm\"uller modular groups,
{\it Lecture Notes in Math.}, No. 1346, Springer-Verlag, 1988, 199-270.

\bibitem[I2]{I2} N. V. Ivanov, {\it Subgroups of Teichm\"uller modular groups,}
Translations of Mathematical Monographs, V. 115,
American Mathematical Society, 1992, xii, 127 pp.

\bibitem[I3]{I3} N. V. Ivanov, Braids and Thurston's classification, Talk at the
Special Session {\it ``Mapping class groups and geometric theory of 
Teichm\"uller spaces''} at the 974-th Meeting of the AMS, Ann Arbor, MI,
March 1-3, 2002.

\bibitem[M]{M} J. D. McCarthy, Normalizers and centralizers 
of pseudo-Anosov mapping classes, Preprint, 1982; revised version 1994.

\end{thebibliography}
\end{document}